\documentclass[reqno]{amsart}
\usepackage{amssymb}
\usepackage{graphicx}
\usepackage{pictex}
\usepackage{mathrsfs}
\usepackage{hyperref}

\begin{document}

\newtheorem{thm}{Theorem}%[section]
\newtheorem{lem}[thm]{Lemma}
\newtheorem{cor}[thm]{Corollary}
\newtheorem{prop}[thm]{Proposition}

\theoremstyle{definition}
\newtheorem{defn}[thm]{Definition}

\theoremstyle{remark}
\newtheorem{rmk}[thm]{Remark}

\newcommand{\nn}{\noindent}

\def\square{\hfill${\vcenter{\vbox{\hrule height.4pt \hbox{\vrule width.4pt
height7pt \kern7pt \vrule width.4pt} \hrule height.4pt}}}$}

\newenvironment{pf}{\noindent {\sc Proof.}\,\,}{\qed \vskip 9pt}

\title{Congruence and Uniqueness of certain Markoff numbers} %
\author{Ying Zhang}
\address{Department of Mathematics, Yangzhou University, Yangzhou 225002, CHINA} %
\email{yingzhang@yzu.edu.cn}
\address{Current Address: IMPA, Estrada Dona Castorina 110, Rio de Janeiro 22460, BRAZIL}
\email{yiing@impa.br}
\date{December 29, 2006}
\thanks{The author is supported by a TWAS-CNPq
postdoctoral fellowship and in part by NSFC grant No. 10671171.}
\begin{abstract}
By making use of only simple facts about congruence, we first show
that every even Markoff number is congruent to $2$ modulo $32$, and
then, generalizing an earlier result of Baragar, establish the
uniqueness for those Markoff numbers $c$ where one of $3c-2$ and
$3c+2$ is a prime power, $4$ times a prime power, or $8$ times a
prime power.
\end{abstract}
\maketitle

\section{Introduction}\label{s:intro}

\nn It is A. A. Markoff who first studied the Diophantine
equation---now known as the {\it Markoff equation}
\begin{eqnarray}\label{eqn:markoff}
a^2+b^2+c^2=3abc
\end{eqnarray}
in late 1870s in his famous work \cite{markoff1880ma} on the minima
of real, indefinite, binary quadratic forms. (For interpretations of
Markoff's work, see \cite{dickson1930book}, \cite{cassels1957book}
and \cite{cusick-flahive}. For its relation with the hyperbolic
geometry of the modular torus, see \cite{cohn1955am} and
\cite{series1985mi}).

The positive integers $a,b,c$ satisfying (\ref{eqn:markoff}) are
particularly important in the work of Markoff, and Frobenius
\cite{frobenius1913} called them the {\it Markoff numbers}.
The solution triples $(a,b,c)$ in positive integers are called the
{\it Markoff triples}. For convenience, we shall not distinguish a
Markoff triple from others obtained by permuting its elements, i.e.,
from its permutation class, and when convenient, usually arrange its
elements in ascending order. We shall call the two Markoff triples
$(1,1,1)$ and $(1,1,2)$ {\it singular}, while all the others {\it
non-singular}. It is easy to show that the elements of a
non-singular Markoff triple are all distinct.

In ascending order of their largest elements, the first 12 Markoff
triples are: $(1,1,1)$, $(1,1,2)$, $(1,2,5)$, $(1,5,13)$,
$(2,5,29)$, $(1,13,34)$, $(1,34,89)$, $(2,29,169)$, $(5,13,194)$,
$(1, 89, 233)$, $(5, 29, 433)$, $(89, 233, 610)$. And the first 40
Markoff numbers are recorded in \cite{sloane}.

That the Markoff equation is particularly interesting lies in the
fact that it is a quadratic equation in each of $a,b$ and $c$, and
hence new solutions can be obtained by a simple process from any
given one, $(a,b,c)$.
To see this, keep $a$ and $b$ fixed and let $c'$ be the other root
of (\ref{eqn:markoff}), regarded as a quadratic equation in $c$.
Since (\ref{eqn:markoff}) can be rewritten as
$c^2-3abc+(a^2+b^2)=0$, we have $c+c'=3ab$ and $cc'=a^2+b^2$. Thus
$c'=3ab-c$ is a positive integer and $(a,b,c')$ is another solution
triple to (\ref{eqn:markoff}) in positive integers, that is, a
Markoff triple. Similarly, we obtain two other Markoff triples
$(a',b,c)$ and $(a,b',c)$. We call these three new Markoff triples
thus obtained the {\it neighbors} of the given one.

In \cite{markoff1880ma}, Markoff demonstrated that every Markoff
triple can be obtained from $(1,1,1)$ by repeatedly generating new
neighbors.

\vskip 6pt

\nn {\sc Theorem} {\bf A} (Markoff \cite{markoff1880ma}).\, {\it
Every Markoff triple can be traced back to $(1,1,1)$ by repeatedly
performing the following operation on Markoff triples:
\begin{eqnarray}\label{eqn:operation}
(a,b,c)\longmapsto(a,b,c'):=(a,b,3ab-c) %
\end{eqnarray}
where the elements of $(a,b,c)$ is arranged so that $a \le b \le c$,
and the elements of $(a,b,c')$ need to be rearranged in ascending
order to perform the next operation.}

\vskip 6pt

Markoff's proof of Theorem A can be found in
\cite[pp.397--398]{markoff1880ma}. A different, simple proof can be
found in \cite[pp.27--28]{cassels1957book}; see also
\cite[pp.17--18]{cusick-flahive}. Another slightly different proof
is given by the author in \cite{zhang2006arxiv}.

The idea of the proof given in \cite{cassels1957book} is that
operation (\ref{eqn:operation}) above reduces the largest elements
of Markoff triples as long as the input triple is non-singular.
Indeed, if $a < b < c$ then
$(c-b)(c'-b)=cc'-(c+c')b+b^2=a^2+2b^2-3ab^2<0$, and hence $c' < b$.
Therefore, after a finite number of steps of reduction, the process
will stop at a singular Markoff triple, which is in fact $(1,1,2)$.
Applying operation (\ref{eqn:operation}) another time then gives
$(1,1,1)$.

As an immediate corollary of Theorem A, we have

\vskip 6pt

\nn {\sc Theorem} {\bf B} (Frobenius \cite{frobenius1913}). \; %

{\it {\rm (a)} The elements of a Markoff triple are pairwise coprime. %

{\rm (b)} Every odd Markoff number is $\equiv 1 \,({\rm mod}\; 4)$.

{\rm (c)} Every even Markoff number is $\equiv 2 \,({\rm mod}\; 8)$.} %

\vskip 6pt

\begin{pf} By (\ref{eqn:markoff}),
$\gcd (a,b)=\gcd (b,c)=\gcd (c,a)=\gcd (a,b,c)$, and by Theorem A,
$\gcd (a,b,c)=\gcd (a,b,3ab-c)=\cdots=\gcd (1,1,1)$ $=1$. This
proves (a).

Since $c(3ab-c)=a^2+b^2$ and $\gcd (a,b)=1$, $c$ is not a multiple
of $4$, and for each prime factor $p$, $-1$ is a quadratic residue
modulo $p$. Since it is well known that $-1$ is not a quadratic
residue modulo a prime $\equiv 3 \,({\rm mod}\; 4)$, each odd prime
factor of $c$ is $\equiv 1 \,({\rm mod}\; 4)$, from which (b) and
(c) follow.
\end{pf}

The following conjecture on the uniqueness of Markoff numbers or
Markoff triples was somewhat hidden in the work of Markoff and was
first mentioned explicitly as a question by G. Frobenius in his 1913
paper \cite{frobenius1913}. It asserts that a Markoff triple is
uniquely determined by its largest element. We shall simply say that
a Markoff number $c$ is {\it unique} if the following holds for $c$.

\vskip 6pt

\nn {\sc The Unicity Conjecture.} \, {\it Suppose \,$(a,b,c)$ and
$(\tilde{a},\tilde{b},c)$ are Markoff triples with $a \le b \le c$
and \,$\tilde{a} \le \tilde{b} \le c$. Then $a=\tilde{a}$ and
$b=\tilde{b}$.}

\vskip 6pt

The conjecture has become widely known when Cassels mentioned it in
\cite[p.33]{cassels1957book}; see also \cite[p.11,
p.26]{cusick-flahive} and \cite[p.188]{conway-guy}. It has been
proved only for some rather special subsets of the Markoff numbers.
The following result for Markoff numbers which are prime powers or
$2$ times prime powers was first proved independently and partly by
Baragar \cite{baragar1996cmb} (for primes and $2$ times primes),
Button \cite{button1998jlms} (for primes but can be easily extended
to prime powers) and Schmutz \cite{schmutz1996ma} (for prime powers
but the proof works also for $2$ times prime powers) using either
algebraic number theory
(\cite{baragar1996cmb},\cite{button1998jlms}) or hyperbolic geometry
(\cite{schmutz1996ma}). A simple, short proof using the hyperbolic
geometry of the modular torus as used by Cohn in \cite{cohn1955am}
has been obtained a bit later but only recently posted by Lang and
Tan \cite{lang-tan2005markoff}. See \cite{zhang2006arxiv} for a
completely elementary proof which uses neither hyperbolic geometry
nor algebraic number theory. A stronger result along the same lines
has been obtained by Button in \cite{button2001jnt}; in particular,
the Markoff numbers which are ``small'' ($\le 10^{35}$) multiples of
prime powers are unique.

\vskip 6pt

\nn {\sc Theorem} {\bf C} (Baragar \cite{baragar1996cmb}; Button
\cite{button1998jlms}; Schmutz \cite{schmutz1996ma}). {\it %
A Markoff number is unique if it is a prime power or $2$ times a
prime power.}

\vskip 6pt

In this paper we first obtain the following simple but sharper
congruence for all even Markoff numbers.

\begin{thm}\label{thm:even}
If $c$ is an even Markoff number then $c \equiv 2 \,({\rm mod}\, 32)$. %
\end{thm}

This congruence is best possible since the first two even Markoff
numbers are $2$ and $34$. And it seems only the congruence $c \equiv
2 \,({\rm mod}\, 8)$ has been previously observed by Frobenius
\cite{frobenius1913}.
As a consequence of Theorem \ref{thm:even}, we see that for an even
Markoff number $c$,
$$3c-2 \equiv 4 \,({\rm mod}\, 32) \quad\quad \text{and} \quad\quad %
3c+2 \equiv 8 \,({\rm mod}\, 32);$$ %
hence $3c-2$ and $3c+2$ are respectively $4$ times and $8$ times an
odd number.

As the other main result of this paper, we then have the following

\begin{thm}\label{thm:unique}
A Markoff number $c$ is unique if one of $3c+2$ and $3c-2$ is a
prime power, $4$ times a prime power, or $8$ times a prime power.
\end{thm}

In the case where one of $3c+2$ and $3c-2$ is a prime or $4$ times a
prime, this has been obtained by Baragar in \cite{baragar1996cmb}
(and earlier by D. Zagier but not published).

The proofs of Theorems \ref{thm:even} and \ref{thm:unique} will be
given in \S \ref{s:pfthm}. Our method is completely self-contained
and elementary in the sense that it uses nothing but very basic
facts on congruence, which we list as Lemmas \ref{lem:odd} and
\ref{lem:roots} in \S \ref{s:lemmas} and include proofs.
It is most important for us to note that Markoff equation
(\ref{eqn:markoff}) can be rewritten as
\begin{eqnarray}
& \hspace{-60pt} (a-b)^2+c^2=ab(3c-2); \label{eqn:3c-2} \\ %
& \hspace{-60pt} (a+b)^2+c^2=ab(3c+2). \label{eqn:3c+2} %
\end{eqnarray}
Actually, the result in Theorem \ref{thm:unique} came to the
author's mind immediately after he saw (\ref{eqn:3c-2}) and
(\ref{eqn:3c+2}) printed in \cite[p.601]{frobenius1913}.

%%%%%%%%%%%%%%%%%%%%%%%%%%%%%%%%%%%%%%%%%%%%%%%%%%%%%%%%%%%%%%%%%%%%%%
\section{Lemmas}\label{s:lemmas}
%%%%%%%%%%%%%%%%%%%%%%%%%%%%%%%%%%%%%%%%%%%%%%%%%%%%%%%%%%%%%%%%%%%%%%

\nn The following two basic facts from elementary number theory will
be used in the proofs of the theorems. Here we include proofs so
that our exposition of this paper be entirely self-contained.

\begin{lem}\label{lem:odd}
If $x$ and $y$ are coprime integers then every odd factor of
$x^2+y^2$ is $\equiv 1 \,({\rm mod}\; 4)$.
\end{lem}

\begin{pf}
It is shown in the proof of Theorem B that every odd prime factor of
$x^2+y^2$ is $\equiv 1 \,({\rm mod}\; 4)$, from which the conclusion
of the lemma follows.
\end{pf}

\begin{lem}\label{lem:roots}
Suppose $m=p^n$ or $2p^n$ for an odd prime $p$ and an integer $n \ge
1$. Then, for any integer $r$ coprime to $m$, the binomial quadratic
equation
\begin{eqnarray}\label{eqn:qbinom}
x^2 + r \equiv 0 \quad ({\rm mod}\; m)
\end{eqnarray}
has at most one integer solution $x$ with $0<x<m/2$.
\end{lem}

\begin{pf}
We prove the lemma for $m=2p^n$ only; the proof for the case where
$m=p^n$ is similar and actually a bit simpler. Suppose
(\ref{eqn:qbinom}) has two integer solutions $x$ and $\tilde{x}$
such that $0<x<\tilde{x}<m/2$. Then $2p^n \mid
(\tilde{x}+x)(\tilde{x}-x)$. Note that $0<\tilde{x}+x<2p^n$ and
$0<\tilde{x}-x<p^n$. If $p \mid \tilde{x}+x$ and $p \mid
\tilde{x}-x$ then $p \mid 2x$; hence $p \mid x$, and consequently $p
\mid r$, a contradiction. Therefore we must have $p^n \mid
\tilde{x}+x$ and $2 \mid \tilde{x}-x$. But then $\tilde{x}+x=p^n$
and $\tilde{x} \equiv x \,({\rm mod}\; 2)$, which implies that $p^n$
is even, again a contradiction. This completes the proof Lemma
\ref{lem:roots}.
\end{pf}

\nn {\bf Remark.}\, One may give a direct proof for Lemma
\ref{lem:roots} using the fact that in this case there is a
primitive root of $m$.

\vskip 6pt

For the proof of Theorem \ref{thm:unique}, we shall also need the
following rough comparison among the elements of a non-singular
Markoff triple.

\begin{lem}\label{lem:twice}
Suppose $(a,b,c) \neq (1,2,5)$ is a Markoff triple with $a<b<c$.
Then
\begin{eqnarray}\label{eqn:twice}
c > 2ab \quad \text{and} \quad b > 2c'a;
\end{eqnarray}
where $c'=3ab-c$; in particular, $c > 2b$ and $b > 2a$.
\end{lem}

\begin{pf}
By Theorem A, every Markoff triple $(a,b,c) \neq (1,2,5)$ can be
obtained by repeatedly generating new neighbors starting from
$(1,2,5)$. Hence we only need to show that if (\ref{eqn:twice})
holds for $(a,b,c)$ then it also holds for the two new neighbors
$(a',b,c)$ and $(a,b',c)$, where $a'=3bc-a$ and $b'=3ca-b$. For this
we only need to check $a' > 2bc$ and $b' > 2ca$, which are very
easy. This proves Lemma \ref{lem:twice}.
\end{pf}

\nn {\bf Remark.}\, It can be seen from the above proof that the
result of Lemma \ref{lem:twice} can be improved, say, as $c > 5ab/2$
and $b > 5c'a/2$ if $(a,b,c) \neq (1,2,5), (2,5,29)$; actually, this
was already known to Frobenius \cite{frobenius1913} with a different
proof. But for our purposes in this paper the weaker result that $c
> 2b$ and $b > 2a$ is enough.

%%%%%%%%%%%%%%%%%%%%%%%%%%%%%%%%%%%%%%%%%%%%%%%%%%%%%%%%%%%%%%%%%%%%%%
\section{Proof of Theorems \ref{thm:even} and \ref{thm:unique}}\label{s:pfthm}
%%%%%%%%%%%%%%%%%%%%%%%%%%%%%%%%%%%%%%%%%%%%%%%%%%%%%%%%%%%%%%%%%%%%%%

\nn {\sc Proof of Theorem \ref{thm:even}.}\,\, Suppose $(a,b,c)$ is
a Markoff triple with $a<b<c$ such that $c$ is even. By Theorem B,
$a \equiv b \equiv 1 \,({\rm mod}\; 4)$ and $c \equiv 2 \,({\rm
mod}\; 8)$, hence $(a-b)/2$ is even and $c/2 \equiv 1 \,({\rm mod}\;
4)$. Since $3c-2 \equiv 4 \,({\rm mod}\; 8)$, we know that
$(3c-2)/4$ is odd. Then (\ref{eqn:3c-2}) gives
\begin{eqnarray}
((b-a)/2)^2+(c/2)^2=ab(3c-2)/4.
\end{eqnarray}
Since $\gcd(c/2,a)=\gcd(c/2,b)=1$ and $\gcd(c/2,(3c-2)/4)=1$, we
know that $c/2$ is coprime with $ab(3c-2)/4$, and consequently,
$(b-a)/2$ and $c/2$ are coprime. Then Lemma \ref{lem:odd} implies
$(3c-2)/4 \equiv 1 \,({\rm mod}\; 4)$, from which follows that $c
\equiv 2 \,({\rm mod}\; 16)$.

Then $3c+2 \equiv 8 \,({\rm mod}\; 16)$ and hence $(3c+2)/8$ is odd.
Now (\ref{eqn:3c+2}) gives
\begin{eqnarray}
((a+b)/2)^2+(c/2)^2=2ab(3c+2)/8.
\end{eqnarray}
Since $\gcd(c/2,a)=\gcd(c/2,b)=1$ and $\gcd(c/2,(3c+2)/4)=1$, we
know that $c/2$ is coprime with $ab(3c+2)/4$, and consequently,
$(a+b)/2$ and $c/2$ are coprime. Then Lemma \ref{lem:odd} implies
$(3c+2)/8 \equiv 1 \,({\rm mod}\; 4)$, from which follows that $c
\equiv 2 \,({\rm mod}\; 32)$. This proves Theorem \ref{thm:even}.
\qed

\vskip 12pt

%%%%%%%%%%%%%%%%%%%%%%%%%%%%%%%%%%%%%%%%%%%%%%%%%%%%%%%%%%%%%%%%%%%%%%%%
%%%%%%%%%%%%%\section{Proof of Theorem \ref{thm:unique}}\label{s:pfthm3}
%%%%%%%%%%%%%%%%%%%%%%%%%%%%%%%%%%%%%%%%%%%%%%%%%%%%%%%%%%%%%%%%%%%%%%%%

\nn \nn {\sc Proof of Theorem \ref{thm:unique}.}\,\, Suppose
\,$(a,b,c)$ and $(\tilde{a},\tilde{b},c)$ are Markoff triples with
$a \le b \le c$ and \,$\tilde{a} \le \tilde{b} \le c$. We proceed to
show that if $3c-2$ or $3c+2$ is of the form $p^n$, $4p^n$ or $8p^n$
for an odd prime $p$ and an integer $n \ge 1$ then $a=\tilde{a}$ and
$b=\tilde{b}$ .

\vskip 6pt

{\sc Case 1.\, $c$ is odd}.

\vskip 6pt

{\sc Subcase 1.1.}\, Suppose $3c-2=p^n$. Write $m=3c-2$. Then
(\ref{eqn:3c-2}) gives
\begin{eqnarray}\label{eqn:1.1}
(b-a)^2+c^2 = abm \equiv 0 \quad ({\rm mod}\; m). %
\end{eqnarray}
Note that $\gcd(c,m)=1$ since $\gcd(c,m) \mid 2$ and $c$ is odd. By
Lemma \ref{lem:twice}
\begin{eqnarray}\label{eqn:1.1m/2}
0< b-a < c/2 -1 < (3c-2)/2 = m/2.
\end{eqnarray}
Since (\ref{eqn:1.1}) and (\ref{eqn:1.1m/2}) are also true for
$(\tilde{a}, \tilde{b}, c)$, Lemma \ref{lem:roots} implies
$b-a=\tilde{b}-\tilde{a}$. Substituting this relation back into
(\ref{eqn:1.1}) and its analog for $(\tilde{a}, \tilde{b}, c)$, one
then obtains $ab=\tilde{a}\tilde{b}$. Hence both $\{-a,b\}$ and
$\{-\tilde{a},\tilde{b}\}$ are the roots of the same quadratic
equation. This implies $a=\tilde{a}$ and $b=\tilde{b}$.

\vskip 6pt

{\sc Subcase 1.2.}\, Suppose $3c+2=p^n$. Write $m=3c+2$. The proof
is similar to that of Subcase 1.1, now using
\begin{eqnarray*}\label{eqn:1.2}
(a+b)^2+c^2 = abm \equiv 0 \quad ({\rm mod}\; m) %
\end{eqnarray*}
and\, $0< a+b < 3c/4 < (3c+2)/2 = m/2$.

\vskip 6pt

{\sc Case 2.\, $c$ is even}.

\vskip 6pt

By Theorem \ref{thm:even}, $3c-2$ is $4$ times an odd and $3c+2$ is
$8$ times an odd. And by Theorem B,\, $a \equiv b \equiv 1 \,({\rm
mod}\; 4)$\, and \,$\tilde{a} \equiv \tilde{b} \equiv 1 \,({\rm
mod}\; 4)$.

\vskip 6pt

{\sc Subcase 2.1.}\, Suppose $3c-2=4p^n$. Write $m=(3c-2)/4=p^n$.
Then (\ref{eqn:3c-2}) gives
\begin{eqnarray}\label{eqn:2.1}
((b-a)/2)^2+(c/2)^2=abm \equiv 0 \quad ({\rm mod}\; m). %
\end{eqnarray}
Since $\gcd(c,3c-2)=2$, we have $\gcd(c/2,m)=1$.  By Lemma
\ref{lem:twice}
\begin{eqnarray}\label{eqn:2.1m/2}
0<(b-a)/2<c/4<(3c-2)/8=m/2.
\end{eqnarray}
Since (\ref{eqn:2.1}) and (\ref{eqn:2.1m/2}) are also true for
$(\tilde{a}, \tilde{b}, c)$, Lemma \ref{lem:roots} implies that
$b-a=\tilde{b}-\tilde{a}$, and consequently,
$ab=\tilde{a}\tilde{b}$. Therefore $a=\tilde{a}$ and $b=\tilde{b}$.

\vskip 6pt

{\sc Subcase 2.2.}\, Suppose $3c+2=8p^n$. Write $m=(3c+2)/4=2p^n$.
The proof is similar to that of Subcase 2.1, now using
\begin{eqnarray*}\label{eqn:2.2}
((a+b)/2)^2+(c/2)^2=abm \equiv 0 \quad ({\rm mod}\; m), %
\end{eqnarray*}
again\, $\gcd(c/2,m)=1$,\, and \,$0<(a+b)/2<3c/8<(3c+2)/8=m/2.$ %

\vskip 6pt

This completes the proof of Theorem \ref{thm:unique}. \qed

\vskip 12pt

\noindent {\bf Acknowledgements.} The author would like to thank Ser
Peow Tan for very helpful conversations and suggestions during his
visit to IMPA, Rio de Janeiro in middle December, 2006.

\end{document}